\documentclass[preprint]{imsart}
\RequirePackage[OT1]{fontenc}
\RequirePackage{amsthm,amsmath}
\RequirePackage[numbers]{natbib}
\RequirePackage[colorlinks,citecolor=blue,urlcolor=blue]{hyperref}
\arxiv{arXiv:math.PR/0000000}
\usepackage{latexsym}
\usepackage{amsmath}
\usepackage{amsthm}
\usepackage{amsfonts}
\usepackage{tikz}
\usetikzlibrary{trees}
\usetikzlibrary{decorations.pathmorphing}
\usetikzlibrary{decorations.markings}

\startlocaldefs
\theoremstyle{plain}\theoremstyle{plain}

\newtheorem{definition}{Definition}

\newtheorem{remark}{Remark}

\newcommand{\field}[1]{\mathbb{#1}}

\newcommand{\R}{\field{R}}
\newcommand{\N}{\field{N}}

\numberwithin{equation}{section}
\numberwithin{lemma}{section}
\numberwithin{theorem}{section}
\numberwithin{corollary}{section}
\numberwithin{remark}{section}
\numberwithin{definition}{section}
\newcommand{\at}{\makeatletter@\makeatother}
\endlocaldefs

\begin{document}
\begin{frontmatter}
\title{The Third and Fourth Moment of the Renormalized Intersection Local Time}
\runtitle{$3^{rd} \& 4^{th}$ Moment}

\begin{aug}
\author{\fnms{Daniel} \snm{H\"of}
\ead[label=e1]{dhoef@gmx.net}}



\affiliation{Provide e.V.}

\address{
2. K\i{}s\i{}m Mah. Anadolu Cad. 16 D: 8\\
Bah\c ce\c sehir\\
TR-34538 \.Istanbul \\
Republic of Turkey\\
\printead{e1}\\
}

\end{aug}

\begin{abstract}
In this article we calculate the third and fourth moment of the renormalized intersection local time of a planar Brownian motion. The third moment is calculated anlaytically, the fourth moment numerically. For the closed planar random walk the third moment of the distribution of the multiple point range is also calculated in leading order.    
\end{abstract}



\end{frontmatter}

\noindent

\noindent
\section{Introduction}\label{intro}
The distribution of the range and the multiple point range of the planar random walk has been an important subject of mathematical research over the last 60 years \citep{Tek}. The first moment of the multiple point range has been calculated by Flatto in \citep{Fla}, the second moment by Jain and Pruitt in \cite{jain}. The leading behaviour (for large length) of the distribution of the appropriately rescaled range of a planar random walk and that of the renormalized intersection local time of the Brownian motion in two dimensions are proportional to each other with a negative real constant of proportionality as established by Le Gall \citep{legall}. This has been extended to a comparable relationship for the multiple point range by Hamana \citep{hamana_ann} \\
Recently I have calculated the characteristic function of the renormalized intersection local time of the planar Brownian motion in terms of certain integrals related to the $\phi^4$ theory in physics \citep{hoef2}. Fortunately these integrals belong to the wider class of Feynman integrals, many of which have been calculated for the famous $g-2$ experiment in high energy physics (for a recent review see \citep{muong}) over the last $30$ years \citep{Schroder,Laporta}. They therefore belong to the best and most studied objects in the field of numerical and analytical calculation of difficult integrals in the history of mankind. \\
In this article I use this wide knowledge collected in decades for the following results: \\
Let $\beta_1$ be the renormalized intersection local time of a planar Brownian motion as defined e.g. in \citep{basschen} and $E(.)$ the expectation value.
Then the second through fourth moments calculated from \citep[eq. 1.1]{hoef2} are:
\begin{equation}\label{eq_m2}
E(\beta_1^2) = \frac{1}{4\pi^2}\cdot \left( 1 + 3 \cdot \zeta_f - \zeta(2) \right) =0.043035\ldots
\end{equation}
\begin{equation}\label{eq_m3}
E(\beta_1^3) = \frac{1}{16 \pi^3} \left( \frac{311 \cdot \zeta(3)}{18} - 4  - 15 \cdot \zeta_f \right) = 0.010178\ldots
\end{equation}
\begin{equation}\label{eq_m4}
E(\beta_1^4) = 0.010063\ldots
\end{equation}
where \eqref{eq_m2} of course reproduces the classical result of Jain and Pruitt \citep{jain}. $\zeta_f$ is given by
\begin{equation}\label{eq_zf}
\zeta_f := \sum_{p=0}^{\infty} \left( \frac{1}{(1 + 3 \cdot p)^2} - \frac{1}{(2 + 3\cdot p)^2} \right) = 0.781302412896486\ldots
\end{equation}
where a high precision calculation of $\zeta_f$ can be done from the formula proven in \citep{Broadhurst}.\\
The numerical values for the skewness and the excess curtosis of the renormalized intersection local time are then:
\begin{equation}\label{eq_g1}
\gamma_1 := \frac{E(\beta_1^3)}{(E(\beta_1^2))^{\frac{3}{2}}} = 1.140051529\ldots
\end{equation}
\begin{equation}\label{eq_g2}
\gamma_2 := \frac{E(\beta_1^4)}{(E(\beta_1^2))^{2}} - 3 =2.4335\ldots
\end{equation}
The asymptotic distribution of $\beta(\omega,\Lambda)$, the multiple point range for the simple closed planar random  walk rescaled appropriately \citep[eq. 8.56]{hoef2}, calculated from \citep[eq. 8.57]{hoef2} has the following second moment
\begin{equation}\label{eq_bo2}
E\left(\beta(\omega,\Lambda)^2 \right) = \frac{1}{8\pi^2}\left(7\zeta(3) - 2\zeta(2) \right) = 0.0649029\ldots
\end{equation}
and third moment 
\begin{equation}\label{eq_bo3}
E\left(\beta(\omega,\Lambda)^3 \right) = -\frac{7 \zeta(3)}{16\pi^3} = -0.016961\ldots
\end{equation}
and therefore the skewness
\begin{equation}\label{eq_ga1c}
\gamma_{1,closed} = -1.0257865\ldots
\end{equation}
The path to these results is straightforward. We proceed from the formulas in \citep[eq. 1.1, eq. 8.57]{hoef2}:  In section \ref{calc} we put together the contributing matrices and coefficients for the moment formulas. In section \ref{integr} we discuss the connection of the integrals to Feynman integrals and generalized gamma functions and give the sources for the integrals involved. In section \ref{mome} we give a short account of the fourth moment and the use of the results for the range and multiple point range of planar random walks. In section \ref{concl} we give a discussion of possible extensions of the results.
\section{Matrices and Coefficients}\label{calc} 
As is clear from \citep[Theorem 1.1 and Theorem 8.4 ]{hoef2} we need to calculate certain coefficients $\mathcal{U}(F)$ and integrals $I(F)$ and $\mathcal{I}(F)$ for all $r \times r$ matrices $F$ which fulfil the following conditions:
\begin{enumerate}
\item $\forall i,j: F_{i,j} \in \N_0$
\item$\forall i: F_{i,i} = 0 $
\item $ \forall j: \sum_{i=0}^r F_{i,j} = \sum_{i=0}^r F_{j,i} = 2$
\item $cof(2\cdot 1_{r \times r}-F) \neq 0$ 
\end{enumerate}
i.e. $F \in \tilde{H}_r(2,\ldots,2)$. \\
In this section we first define an equivalence relation on $\tilde{H}_r(2,\ldots,2)$ to reduce the number of matrices $F$ which enter the calculations and then give an overview of all relevant matrices for $r=2,3,4$ together with their coefficients $\mathcal{U}(F)$. \\ 
For any matrix $F \in \tilde{H}_r(2,\ldots,2)$ and $\sigma \in S_r$ we can define the matrix $F^{\sigma}$ by 
\begin{equation}\label{eq_fsigma}
F^{\sigma}_{i,j} := F_{\sigma(i),\sigma(j)}
\end{equation}
Then it is immediately clear that $F^{\sigma} \in \tilde{H}_r(2,\ldots,2)$ and $\mathcal{U}(F^{\sigma}) = \mathcal{U}(F)$ and $I(F^{\sigma}) = I(F)$ and $\mathcal{I}(F^{\sigma}) = \mathcal{I}(F)$.\\
For $F \in \tilde{H}_r(2,\ldots,2)$ we define the normal subgroup 
$
S_F \subset S_r$ by 
\begin{equation}
\sigma \in S_F \iff F^{\sigma} = F
\end{equation}
On $\tilde{H}_r(2,\ldots,2)$ we can define an equivalence relation $\sim$ by 
\begin{equation}
F \sim F^{'} \iff \exists \sigma \in S_r: F' = F^{\sigma}
\end{equation}
We then have 
\begin{equation}
g(F) := \frac{r!}{\#(S_F)}
\end{equation}
different representatives in each class (where $\#$ denotes the cardinality of a set). 
We then define the set of equivalence classes of $\tilde{H}_r(2,\ldots,2)$ by $\bar{H}_r$. Summing over $\tilde{H}_r(2,\ldots,2)$ can then be replaced by summing over $\bar{H}_r$ with the weights $g(F)$. It is of course quite easy to find all $F \in \bar{H}_r$ for a given $r$ by hand (for small $r$) and by simple computer programs for bigger $r$. \\
We start with $r=2$. There is just one matrix $f_1 \in \bar{H}_2$:  
\begin{equation}
f_1 := \left(
\begin{array}[2]{cc}
	0 & 2 \\
	2 & 0 \\
\end{array} \right)
\end{equation}
and it has $g(f_1) = 1$.\\
For $r=3$ there are two matrices in $\bar{H}_3$:
\begin{equation}
f_2 := \left(
\begin{array}[3]{ccc}
	0 & 1 & 1 \\
	1 & 0 & 1 \\
	1 & 1 & 0 \\
\end{array} \right)
f_3 := \left(
\begin{array}[3]{ccc}
	0 & 2 & 0 \\
	0 & 0 & 2 \\
	2& 0 & 0 \\
\end{array} \right)
\end{equation}
with $g(f_2) = 1$ 
and  $g(f_3) = 2$.\\
For $r=4 \Rightarrow \#\bar{H}_4 = 5$, the matrices are 
\begin{equation}
f_4 := \left(
\begin{array}[4]{cccc}
	0 & 0 & 0 & 2\\
	0 & 0 & 2 & 0\\
	1& 1 & 0 & 0\\
	1& 1 & 0 & 0 \\
\end{array} \right)
f_5 := \left(
\begin{array}[4]{cccc}
	0 & 0 & 0& 2\\
	0 & 0 & 2& 0\\
	2& 0 & 0 & 0\\
	0& 2 & 0 & 0 \\
\end{array} \right)
f_6 := \left(
\begin{array}[4]{cccc}
	0 & 0 & 0& 2\\
	1& 0 & 1 & 0\\
	1& 1 & 0 & 0\\
	0& 1 & 1 & 0 \\
\end{array} \right)
\end{equation}


\begin{equation}
f_7 := \left(
\begin{array}[4]{cccc}
	0 & 0 & 1& 1\\
	0 & 0 & 1 & 1\\
	1& 1 & 0 & 0\\
	1& 1 & 0 & 0 \\
\end{array} \right)
f_8 := \left(
\begin{array}[4]{cccc}
	0 & 0 & 1& 1\\
	1& 0 & 0 & 1\\
	1& 1 & 0 & 0\\
	0& 1 & 1 & 0 \\
\end{array} \right)
\end{equation}


The following table now gives the coefficients as defined in \citep[Theorem 1.1 and Theorem 8.4 ]{hoef2}.\\
\\
\renewcommand{\arraystretch}{1.5}
\begin{tabular}[h]{c|c|c|c|c}
$F$ & $g(F)$ & $ cof(2\cdot 1_{r \times r} - F)$ & $M(F)$ & $g(F) \cdot \mathcal{U}(F)$  \\
\hline 
$f_1$ & 1 & 2 & 4 & $\frac{1}{2} $ \\
\hline 
$f_2$ & 1 & 3 & 1 & 3 \\
\hline 
$f_3$ & 2 & 4 & 8 & 1 \\
\hline 
$f_4$ & 12 & 4 & 4 & 12\\
\hline 
$f_5$ & 6 & 8 & 16 & 3 \\
\hline 
$f_6$ & 12 & 6 & 2 & 36 \\
\hline 
$f_7$ & 3 & 4 & 1 & 12\\
\hline 
$f_8$ & 6 & 5 & 1 & 30\\
\end{tabular}
\\
\\
\section{Integrals and Graphs}\label{integr}
In this section we first discuss the class of Feynman integrals relevant for our problem. We then give the relation between them and the integrals $I(F)$ and $\mathcal{I}(F)$. We then proceed with the calculation of all integrals needed for the calculation of the second, third and fourth moment.\\ 
For a general introduction into Feynman integrals see e.g. \citep[ch. 6 - 8]{itzykson}.\\
Any $F \in \tilde{H}_r(2,\ldots,2)$ interpreted as adjacency matrix defines a damfree Eulerian multigraph $G(F)$. \\Now for any connected damfree finite multigraph $G$ the (fully massive bosonic) Feynman integral  $I_G$ in $d$ dimensions with external moments zero can be defined as \citep[eq. 8.5]{itzykson}
\begin{equation}
I_G(d) := \frac{1}{(4\cdot \pi)^{L(G)\cdot \frac{d}{2}}} \Gamma_G\left(1-\frac{d}{2}\right)
\end{equation}
where  
\begin{equation}
L(G) := \#E(G) - \#V(G) + 1
\end{equation}
is the so called number of loops of $G$ and $E(G)$ and $V(G)$ the set of edges and vertices of $G$ respectively and $\#$ denotes the number of elements. The generalized gamma function $\Gamma_G$ of a connected damfree multigraph $G$ for $\Re(n) \geq 0$ is defined by
\begin{equation}
\Gamma_G(n) := \int \prod_{e \in E(G)} d\alpha_e \exp\left( - \sum_{e \in E(G)} \alpha_e \right) \cdot \mathcal{P}_G(\alpha)^{n-1}
\end{equation}
where the integration over the real variables $\alpha_e$ for each edge $e \in E(G)$ is performed over the interval $[0,\infty]$ and $\mathcal{P}_G(\alpha)$ is the so called Kirchhoff-Symanzik polynomial defined by
\begin{equation}\label{eq_poly}
\mathcal{P}_G(\alpha) := \sum_{T \in sp(G)} \left( \prod_{e \in E(G \setminus T)} \alpha_e \right)
\end{equation}
where $sp(G)$ is the set of spanning trees of $G$. \\
\begin{definition}\label{de_ge}
We call two damfree connected graphs $G_1$ and $G_2$ gamma-equivalent
if there are bijective maps $h_1 : E(G_1) \mapsto E(G_2)$ and $h_2 :sp(G_1) \mapsto sp(G_2)$ such that for any spanning tree $T \in sp(G_1)$  
\begin{equation}
e \notin T \Leftrightarrow h_1(e) \notin h_2(T)
\end{equation}
\end{definition}
\begin{remark}
From equation \eqref{eq_poly} we immediately see that for two gamma-equivalent graphs $G_1$ and $G_2$  
\begin{equation}
\Gamma_{G_1} = \Gamma_{G_2}
\end{equation}
\end{remark}
For many generalized gamma functions functional relations have been found from integration by parts. 
For example for $f_1$ and $f_2$ we find from the literature \citep[p. 12]{Wegner_pp}
\begin{equation}\label{eq_gf1}
\Gamma_{G(f_1)} (n + 1) =  \frac{3(3n + 2)(3n + 1)(2n + 1)n}{32(n + 1)} \Gamma_{G(f_1)}(n)  + \frac{11 n + 8}{8 (n+ 1)} \Gamma(n + 1)^3
\end{equation}
and from  \citep[p. 17]{Wegner_pp}
\begin{multline}\label{eq_gf2}
\Gamma_{G(f_2)}(n + 1) = \frac{1}{(105 + 286n + 252n^2+ 72n^3)} \cdot \Biggl( \\
 \frac{ 4n \cdot (15 + 137n + 510 n^2 + 988n^3 + 1048n^4 + 576n^5 + 128n^6)}{3} \Gamma_{G(f_2)}(n) + \\
\frac{n(810 + 6905n + 22363n^2 + 34450n^3 + 25268n^4 + 7080n^5)}{32} \Gamma(n + 1) \Gamma_{G(f_1)}(n) \\
+ \frac{840 + 2893n + 3228n^2 + 1172n^3}{8} \Gamma(n+1)^4  \Biggr)
\end{multline} \\
By virtue of functional relations as the above ones the generalized gamma functions can be extended to the full complex plane up to poles. For integers $d >2$ the Laurent expansions for the Feynman integrals around integer $d$ have been calculated in many cases \citep{Schroder, Laporta}.\\
Now how are $I(F)$ and $\mathcal{I}(F)$ related to the above Feynman integrals?
Because of \citep[eq. 7.24]{kleinert}
\begin{equation}\label{eq_prop}
\frac{K_0(m \left| x \right|)}{2\pi} = \int_{\R^2} \frac{e^{iq\cdot x}}{q^2 + m^2} d^2q
\end{equation}
for $x \in \R^2 \setminus \{0 \}$ and real strictly positive $m \in (0,\infty)$. Therefore the integrals $I(F)$ can be related to the generalized gamma functions $\Gamma_{G(F)}(n)$ according to \cite[eq. 8.4, 8.5]{itzykson} by 
\begin{equation}
I(F) = \frac{4^{\#E(G)}}{(4\cdot \pi)^{L(G)}} \Gamma_{G(F)}\left(0\right)
\end{equation}
where of course $G = G(F)$. 
By the same token we find 
\begin{equation}
\mathcal{I}(F) = \frac{4^{\#E(G) - 1}}{(4\cdot \pi)^{L(G)-1}} \left( \sum_{e \in E(G)} \Gamma_{G(F) \setminus e} (0)  \right)
\end{equation} 
where again $G = G(F)$ and $G(F) \setminus e$ denotes the graph which one receives by removing the edge $e$ from $G(F)$. As $G(F)$ is Eulerian $G(F) \setminus e$ is connected. \\
We can now discuss the results for $f_1,\ldots,f_8$. 
To reduce the number of calculations let us first note that for a given matrix $F$ the integral $I(F)$ and the sum of the integrals $\mathcal{I}(F)$ only depends on the sums $F_{i,j} + F_{j,i}$. We will therefore call matrices $F_1$ and $F_2$ transpose equivalent  if there is a $\sigma \in S_r$ such that $F_1^{\sigma} + (F_1^{\sigma})^t = F_2 + F_2^t$. ($F_1^{\sigma}$ as defined in \eqref{eq_fsigma} and $t$ the sign for transposition). We note that if $F_1$ and $F_2$ are transpose equivalent then $I(F_1) = I(F_2)$ and $\mathcal{I}(F_1) = \mathcal{I}(F_2)$ . We see that $f_2$ and $f_3$ are transpose equivalent. $f_5$ and $f_7$ are transpose equivalent too and so are  $f_6$ and $f_8$ . \\
From \eqref{eq_gf1} together with \citep[eq. 6.13]{Schroder} we find 
\begin{equation}
I(f_1) = \frac{28}{\pi^3} \cdot  \zeta(3)
\end{equation}
From \citep[p.6]{Bailey} we find 
\begin{equation}
\mathcal{I}(f_1) = 4 \cdot \left( \frac{2}{\pi} \right)^3 \cdot 2 \pi \int_0^{\infty} x \cdot K_0(x)^3 dx = 4\cdot \left( \frac{2}{\pi}\right)^3 \cdot \left( \frac{3\pi}{2} \zeta_f \right) = \frac{48}{\pi^2} \cdot \zeta_f
\end{equation}
where $\zeta_f$ was defined in \eqref{eq_zf}.
From \eqref{eq_gf2} together with \citep[eq. 6.38]{Schroder} we find 
\begin{equation}
I(f_2) = I(f_3) = \frac{48}{\pi^4}\cdot \zeta(3)
\end{equation}
From \citep{Wegner} we also find 
\begin{equation}
\mathcal{I}(f_2) = \mathcal{I}(f_3) = 6 \cdot \frac{4^2}{\pi^3}\cdot \left(\frac{11}{9} \zeta(3) \right) = \frac{352}{3\pi^3} \cdot \zeta(3)
\end{equation}
The graph $G(f_4)$ has two edges which have no parallel and two triples of edges which are parallel to each other. We can therefore write 
\begin{equation}
\mathcal{I}(f_4) = \left( \frac{2}{\pi} \right)^7 \left( 2 \cdot T_U + 6\cdot T_D \right)
\end{equation}
where the contribution $T_U$ from removing an edge which has no parallel edge can now be calculated 
like the contribution to $\mathcal{I}(f_1)$ resulting in 
\begin{equation}
T_U = 2\pi \cdot \left(\frac{3\pi}{2} \zeta_f \right)^2 =\frac{9  \pi^3}{2} \cdot \zeta_f^2
\end{equation} 
For the contribution $T_D$ from removing an edge which has two parallel edges in the graph $G(f_4)$ we first note that the graph belonging to $T_D$ is gamma-equivalent to a graph $\hat{V}_3$ which results from $V_3$ in \citep[Fig. 1, p. ii]{Laporta} by replacing the edge with no paralell in $V_3$ with two edges connected with a new vertex.
We get
\begin{equation}\label{eq_td}
T_D := \int_{\R^2}d^2x \int_{\R^2} d^2y \left(
K_0(\left| x \right|)^3 \cdot  K_0(\left| y \right|)^2 \cdot \pi \cdot  \left| x - y \right|  \cdot  K_1(\left| x - y\right|) \right)
\end{equation}
 where we have used that for $x \in \R^2 \setminus\{0\}$ the equation 
\begin{equation}\label{eq_k0k1}
\int_{\R^2} d^2y \cdot K_0(\left| y \right|) \cdot K_0(\left| x - y \right|) = \pi \left| x \right| \cdot K_1 (\left| x \right|)
\end{equation}
is true. \eqref{eq_k0k1} can be derived from \eqref{eq_prop}  by differentiation in $m^2$ and setting $m = 1$ afterwards. \\
For $\mathcal{I}(f_5)$ we find 
\begin{equation}
\mathcal{I}(f_5) = \frac{64}{\pi^4} \cdot \left( 8 \Gamma_{V_8}(0)\right)
\end{equation}
where $V_8$ is the graph with that name in  \citep[Fig. 1, p. ii]{Laporta}.
For $\mathcal{I}(f_6)$ we find 
\begin{equation}
\mathcal{I}(f_6) = \frac{64}{\pi^4} \cdot \left( 4 \Gamma_{V_7}(0) + 4 \Gamma_{V_5}(0)\right)
\end{equation}
Again the names for the graphs are taken from \citep[Fig. 1, p. ii]{Laporta}. $V_7$ is gamma-equivalent (as defined in Definition \ref{de_ge}) to the graph $G(f_6) \setminus e$ where $e$ is any of the  4 edges which has no paralell in $G(f_6)$. \\
The table below summarizes the integrals of the graphs: \\
\begin{tabular}[h]{c|c|c}
$F$ & $I(F)$ &  $\mathcal{I}(F)$  \\
\hline 
$f_1$ & $\frac{28}{\pi^3}\cdot \zeta(3)$ & $\frac{48}{\pi^2} \cdot \zeta_f $ \\
\hline 
$f_2$ & $\frac{48}{\pi^4} \cdot \zeta(3)$ & $\frac{352}{3 \pi^3} \cdot \zeta(3) $\\
\hline 
$f_3$ & $I(f_2)$&  $\mathcal{I}(f_2)$ \\
\hline 
$f_4$ &  &  $\left(\frac{2}{\pi}\right)^7 \cdot  \left( 9 \cdot \pi^3 \cdot \zeta_f^2 + 6 \cdot T_D \right) $ \\
\hline 
$f_5$ &  & $\frac{64}{\pi^4} \cdot (8 \cdot \Gamma_{V_8}(0))$ \\
\hline 
$f_6$ &  & $\frac{64}{\pi^4} \cdot (4 \cdot \Gamma_{V_7}(0) + 4 \cdot \Gamma_{V_5}(0))$ \\
\hline 
$f_7$ &  & $\mathcal{I}(f_5)$ \\
\hline 
$f_8$  & & $\mathcal{I}(f_6)$ \\
\end{tabular}
\\
\\
The values of $\Gamma_{V_i}(0)$  for $i = 5,7,8$  (and all other $i$) have been calculated by Laporta to the precision of $1200$ digits \citep{Laporta2}. $T_D$ can be calculated numerically from \eqref{eq_td} after reducing it to the three-dimensional integral.  
\begin{multline}\label{eq_tdnum}
T_D = 2\pi \int_{0}^{\infty} dr_1 \int_{0}^{\infty} dr_2 \int_{0}^{2\pi}d\phi \Biggl(
r_1 \cdot K_0(r_1)^3 \cdot r_2 \cdot  K_0(r_2)^2 \cdot \\ \pi \cdot  \sqrt{r_1^2 + r_2^2 - 2r_1r_2cos(\phi)}  \cdot  K_1\left(\sqrt{r_1^2 + r_2^2 - 2r_1r_2cos(\phi)}\right)\Biggr)
\end{multline}
$T_D$ could also be calculated from recurrence relations for $\Gamma_{V_3}$ to high precision. This is due to the fact that the graph belonging to $T_D$ is gamma-equivalent to a graph $\hat{V}_3$ which results from $V_3$ by replacing the edge with no paralell in $V_3$ with two edges connected with a new vertex. There are functional relations between $\Gamma_{\hat{V}_3} $ and $\Gamma_{V_3}$ which are at the heart of the high precision numerical algorithm of Laporta as described in \citep{Laporta_alg} \\
$\Gamma_{V_7}(0)$ and $\Gamma_{V_8}(0)$ also can be calculated from the recurrence relations in \citep[p.17-18]{Wegner_pp} together with \citep[eq. 7, eq. 8]{Laporta}.  
\section{The moments}\label{mome}
From the above the moments as given in equations \eqref{eq_m2}, \eqref{eq_m3}, \eqref{eq_bo2} and \eqref{eq_bo3} follow by simple computer algebra, for the fourth moment we find:
\begin{multline}
E(\beta_1^4) = \frac{1}{16 \pi^4}\cdot \Biggl( 11\cdot (\Gamma_{V_5}(0) + \Gamma_{V_7}(0)) + 5\cdot \Gamma_{V_8}(0)  + \frac{6}{ \pi^3} \cdot T_D \\+ 9 - \pi^2  + \frac{\pi^4}{60} + (37  - 3\cdot \pi^2) \cdot \zeta_f  + 9 \cdot \zeta_f^2 -  \frac{1243 \cdot \zeta(3)}{54}  \Biggr) = 0.010063\ldots
\end{multline}
For the numerical evaluation of the fourth moment we have calculated $T_D$ by numerical integration and used the values of Laporta \citep{Laporta2} for $\Gamma_{V_i}(0)$ for $i = 5, 7, 8$.\\
From the above we also know for strongly aperiodic planar random walks $w$ with individual steps of mean zero and finite covariance matrix $\Xi^2$  (see \citep[p. 601]{hamana_ann}): the leading term of the distribution of the range $R_{n}$ as defined in \citep[p. 31]{jain} for large length $n$ rescaled as 
\begin{equation}
\frac{\ln(n)^2}{n} \cdot \left( R_{n} - E(R_{n}) \right)
\end{equation}
and of the multiple point range $Q_n^{(p)}$ as defined in \citep[p. 598]{hamana_ann} rescaled as
\begin{equation}
\frac{\ln(n)^3}{n} \cdot \left( Q_{n}^{(p)} - E(Q_{n}^{(p)}) \right)
\end{equation}
have a skewness of  $-\gamma_1$  (as given in \eqref{eq_g1}) and excess curtosis of  $\gamma_2$ (as given in \eqref{eq_g2}). This is true as for large length they converge to the distribution of the intersection local time times  the negative constants $-4 \cdot \pi^2 \cdot \det(\Xi)$ \citep[Theorem 6.1]{legall} for the range and $-16 \cdot \pi^3 \cdot \det(\Xi^2)$ \citep[Theorem 3.5]{hamana_ann} for the mutliple point range respectively.
\section{Conclusion}\label{concl}
We have calculated the third and fourth moment of the renormalized intersection local time of the planar Brownian motion hitherto unknown. It gives us a broad understanding of its shape. The fact that $\gamma_1 > 0$ and $\gamma_2 > 0$ is in line with the fact, that the distribution of $\beta_1$ has an exponential tail for large positive $\beta_1$ and a double exponential tail for large negative $\beta_1$ as proven in \citep[Theorem 1.2 and Theorem 1.3]{basschen}. For illustration purposes: an example of a distribution which has a similar skewness and excess curtosis and a broadly similar behaviour in the tails is the Gumbel distribution.\\It is also interesting to see that the value of the skewness of the asymptotic distribution of the multiple point range for the unrestricted random walk is only slightly different from its counterpart for the closed random walk.\\
This work could be extended in two directions:
\begin{enumerate}
\item  higher moments could be calculated to get an even better understanding of the distribution, especially after an asymptotic calculation of the moments for large $r$ has been done. 
\item higher corrections for large length to the second, third and fourth moments for the distribution of the multiple point range could be calculated. (see e.g. \citep[eq. 8.20, eq. 9.46]{hoef2} for the first correction to the second moment). This would make it possible to get a good comparison between the calculated moments of the distribution and real life planar random walks.
\end{enumerate}
 The relevant matrices can be obtained from a simple computer program in both cases. For a similar problem of $\phi^4$ theory they have been obtained up to order $7$ \citep[ch. 14.2, Fig. 14.3]{kleinert} which is equivalent to the $7^{th}$ moment in our case. The integrals on the other hand are more difficult to calculate analytically or numerically. For higher moments a Monte-Carlo integration of generalized gamma functions for $n=0$ seems a possibility. For corrections to the second, third and fourth moment for large length numerical calculation of the integrals can be done by standard methods.  But it may also be that the next decade will see a major breakthrough in the analytical calculation of the generalized gamma functions. They do play a fundamental role in physics but as we have seen here also for the distribution of the intersection local time of the planar Brownian motion. \\
I am very grateful for the generous sharing of the numerical calculations of $\Gamma_{V_i}(0)$  for $i=5,7,8$ by Prof. S. Laporta without which this work could not have been completed.   
\bibliographystyle{imsart-number}
\bibliography{hoef_34}

\begin{thebibliography}{18}

\bibitem{Bailey}
\begin{barticle}[author]
\bauthor{\bsnm{Bailey},~\bfnm{David~H.}\binits{D.~H.}},
  \bauthor{\bsnm{Borwein},~\bfnm{Jonathan~M.}\binits{J.~M.}},
  \bauthor{\bsnm{Broadhurst},~\bfnm{David}\binits{D.}} \AND
  \bauthor{\bsnm{Glasser},~\bfnm{M.~L.}\binits{M.~L.}}
(\byear{2008}).
\btitle{Elliptic integral evaluations of {B}essel moments and applications}.
\bjournal{J. Phys. A}
\bvolume{41}
\bpages{205203, 46}.
\bdoi{10.1088/1751-8113/41/20/205203}
\bmrnumber{2450513 (2010c:33011)}
\end{barticle}
\endbibitem

\bibitem{basschen}
\begin{barticle}[author]
\bauthor{\bsnm{Bass},~\bfnm{Richard~F.}\binits{R.~F.}} \AND
  \bauthor{\bsnm{Chen},~\bfnm{Xia}\binits{X.}}
(\byear{2004}).
\btitle{Self-intersection local time: critical exponent, large deviations, and
  laws of the iterated logarithm}.
\bjournal{Ann. Probab.}
\bvolume{32}
\bpages{3221--3247}.
\bdoi{10.1214/009117904000000504}
\bmrnumber{2094444 (2005i:60149)}
\end{barticle}
\endbibitem

\bibitem{Broadhurst}
\begin{barticle}[author]
\bauthor{\bsnm{Broadhurst},~\bfnm{David~J.}\binits{D.~J.}}
(\byear{1999}).
\btitle{{Massive three - loop Feynman diagrams reducible to SC* primitives of
  algebras of the sixth root of unity}}.
\bjournal{Eur.Phys.J.}
\bvolume{C8}
\bpages{311-333}.
\bdoi{10.1007/s100529900935}
\end{barticle}
\endbibitem

\bibitem{Tek}
\begin{barticle}[author]
\bauthor{\bsnm{Dvoretzky},~\bfnm{A.}\binits{A.}},
  \bauthor{\bsnm{Erd{\"o}s},~\bfnm{P.}\binits{P.}} \AND
  \bauthor{\bsnm{Kakutani},~\bfnm{S.}\binits{S.}}
(\byear{1950}).
\btitle{Double points of paths of {B}rownian motion in {$n$}-space}.
\bjournal{Acta Sci. Math. Szeged}
\bvolume{12}
\bpages{75--81}.
\bmrnumber{0034972 (11,671e)}
\end{barticle}
\endbibitem

\bibitem{Fla}
\begin{barticle}[author]
\bauthor{\bsnm{Flatto},~\bfnm{Leopold}\binits{L.}}
(\byear{1976}).
\btitle{The multiple range of two-dimensional recurrent walk}.
\bjournal{Ann. Probability}
\bvolume{4}
\bpages{229--248}.
\bmrnumber{0431388 (55 \#\#4388)}
\end{barticle}
\endbibitem

\bibitem{hamana_ann}
\begin{barticle}[author]
\bauthor{\bsnm{Hamana},~\bfnm{Yuji}\binits{Y.}}
(\byear{1997}).
\btitle{The fluctuation result for the multiple point range of two-dimensional
  recurrent random walks}.
\bjournal{Ann. Probab.}
\bvolume{25}
\bpages{598--639}.
\bdoi{10.1214/aop/1024404413}
\bmrnumber{1434120 (98f:60136)}
\end{barticle}
\endbibitem

\bibitem{hoef2}
\begin{barticle}[author]
\bauthor{\bsnm{{Hoef}},~\bfnm{D.}\binits{D.}}
(\byear{2013}).
\btitle{{The Characteristic Function of the Renormalized Intersection Local
  Time of the Planar Brownian Motion}}.
\bjournal{ArXiv e-prints}.
\end{barticle}
\endbibitem

\bibitem{itzykson}
\begin{bbook}[author]
\bauthor{\bsnm{Itzykson},~\bfnm{C.}\binits{C.}} \AND
  \bauthor{\bsnm{Zuber},~\bfnm{J.~B.}\binits{J.~B.}}
(\byear{1988}).
\btitle{Quantum Field Theory}.
\bpublisher{Mc Graw Hill}, \baddress{New York}.
\end{bbook}
\endbibitem

\bibitem{jain}
\begin{binproceedings}[author]
\bauthor{\bsnm{Jain},~\bfnm{N.~C.}\binits{N.~C.}} \AND
  \bauthor{\bsnm{Pruitt},~\bfnm{W.~E.}\binits{W.~E.}}
(\byear{1973}).
\btitle{The range of random walks}.
In \bbooktitle{Proc. Sixth Berkeley Symp. Math. Stat. Probab.}
\bpages{31--50}.
\end{binproceedings}
\endbibitem

\bibitem{muong}
\begin{barticle}[author]
\bauthor{\bsnm{Jegerlehner},~\bfnm{Fred}\binits{F.}} \AND
  \bauthor{\bsnm{Nyffeler},~\bfnm{Andreas}\binits{A.}}
(\byear{2009}).
\btitle{{The Muon g-2}}.
\bjournal{Phys.Rept.}
\bvolume{477}
\bpages{1-110}.
\bdoi{10.1016/j.physrep.2009.04.003}
\end{barticle}
\endbibitem

\bibitem{kleinert}
\begin{bbook}[author]
\bauthor{\bsnm{Kleinert},~\bfnm{H.}\binits{H.}} \AND
  \bauthor{\bsnm{Schulte-Frohlinde},~\bfnm{V.}\binits{V.}}
(\byear{2001}).
\btitle{Critical Properties of $\Phi^4$-theories}.
\bpublisher{World Scientific}.
\end{bbook}
\endbibitem

\bibitem{Laporta_alg}
\begin{barticle}[author]
\bauthor{\bsnm{Laporta},~\bfnm{S.}\binits{S.}}
(\byear{2000}).
\btitle{High-precision calculation of multiloop {F}eynman integrals by
  difference equations}.
\bjournal{Internat. J. Modern Phys. A}
\bvolume{15}
\bpages{5087--5159}.
\bdoi{10.1142/S0217751X00002157}
\bmrnumber{1811366 (2001m:81194)}
\end{barticle}
\endbibitem

\bibitem{Laporta}
\begin{barticle}[author]
\bauthor{\bsnm{Laporta},~\bfnm{S.}\binits{S.}}
(\byear{2002}).
\btitle{High-precision $\epsilon$-expansions of massive four-loop vacuum
  bubbles}.
\bjournal{Physics Letters B}
\bvolume{549}
\bpages{115--122}.
\bdoi{10.1016/s0370-2693(02)02910-6}
\end{barticle}
\endbibitem

\bibitem{Laporta2}
\begin{bunpublished}[author]
\bauthor{\bsnm{Laporta},~\bfnm{S.}\binits{S.}}
(\byear{2014}).
\btitle{private e-mail on the results of the high precision numerical
  calculation of massive four-loop vacuum bubbles in two dimensions}.
\end{bunpublished}
\endbibitem

\bibitem{legall}
\begin{barticle}[author]
\bauthor{\bsnm{Le~Gall},~\bfnm{J.~F.}\binits{J.~F.}}
(\byear{1986}).
\btitle{Propri\'et\'es d'intersection des marches al\'eatoires. {I}.
  {C}onvergence vers le temps local d'intersection}.
\bjournal{Comm. Math. Phys.}
\bvolume{104}
\bpages{471--507}.
\bmrnumber{840748 (88d:60182)}
\end{barticle}
\endbibitem

\bibitem{Schroder}
\begin{barticle}[author]
\bauthor{\bsnm{Schroder},~\bfnm{Y.}\binits{Y.}} \AND
  \bauthor{\bsnm{Vuorinen},~\bfnm{A.}\binits{A.}}
(\byear{2005}).
\btitle{{High-precision epsilon expansions of single-mass-scale four-loop
  vacuum bubbles}}.
\bjournal{JHEP}
\bvolume{0506}
\bpages{051}.
\bdoi{10.1088/1126-6708/2005/06/051}
\end{barticle}
\endbibitem

\bibitem{Wegner}
\begin{barticle}[author]
\bauthor{\bsnm{Wegner},~\bfnm{Franz}\binits{F.}}
(\byear{1987}).
\btitle{Anomalous dimensions for the nonlinear sigma-model in {$2+\varepsilon$}
  dimensions. {I}, {II}}.
\bjournal{Nuclear Phys. B}
\bvolume{280}
\bpages{193--209, 210--224}.
\bdoi{10.1016/0550-3213(87)90144-1}
\bmrnumber{878909 (88i:81124)}
\end{barticle}
\endbibitem

\bibitem{Wegner_pp}
\begin{barticle}[author]
\bauthor{\bsnm{Wegner},~\bfnm{Franz}\binits{F.}}
(\byear{1988}).
\btitle{Four-Loop Order $\beta$-Function of Nonlinear $\sigma$-Models in
  Symmetric Spaces}.
\bjournal{Preprint des Sonderforschungsbereichs 123}
\bvolume{476}.
\end{barticle}
\endbibitem

\end{thebibliography}
\end{document}